\documentclass{amsart}
\usepackage{setspace}

\usepackage{amsmath,amsthm,amsfonts,amscd,amssymb}

\newtheorem{thm}{Theorem}[section]
\newtheorem{cor}[thm]{Corollary}
\newtheorem{lem}[thm]{Lemma}

\theoremstyle{definition}

\theoremstyle{remark}

\begin{document}

\title{Small zeros of quadratic forms over $\overline{\mathbb Q}$}
\author{Lenny Fukshansky}

\address{Department of Mathematics, Mailstop 3368, Texas A\&M University, College Station, Texas 77843-3368}
\email{lenny@math.tamu.edu}

\begin{abstract}
Let $N \geq 2$ be an integer, $F$ a quadratic form in $N$ variables over $\overline{\mathbb Q}$, and $Z \subseteq \overline{\mathbb Q}^N$ an $L$-dimensional subspace, $1 \leq L \leq N$. We prove the existence of a small-height maximal totally isotropic subspace of the bilinear space $(Z,F)$. This provides an analogue over $\overline{\mathbb Q}$ of a well-known theorem of Vaaler proved over number fields. We use our result to prove an effective version of Witt decomposition for a bilinear space over $\overline{\mathbb Q}$. We also include some related effective results on orthogonal decomposition and structure of isometries for a bilinear space over $\overline{\mathbb Q}$. This extends previous results of the author over number fields. All bounds on height are explicit. 
\end{abstract}

\maketitle

\def\A{{\mathcal A}}
\def\B{{\mathcal B}}
\def\C{{\mathcal C}}
\def\D{{\mathcal D}}
\def\F{{\mathcal F}}
\def\x{{\mathcal H}}
\def\I{{\mathcal I}}
\def\J{{\mathcal J}}
\def\K{{\mathcal K}}
\def\L{{\mathcal L}}
\def\M{{\mathcal M}}
\def\R{{\mathcal R}}
\def\s{{\mathcal S}}
\def\V{{\mathcal V}}
\def\X{{\mathcal X}}
\def\Y{{\mathcal Y}}
\def\Z{{\mathcal Z}}
\def\H{{\mathcal H}}
\def\OO{{\mathcal O}}
\def\cee{{\mathbb C}}
\def\Nn{{\mathbb N}}
\def\pee{{\mathbb P}}
\def\que{{\mathbb Q}}
\def\real{{\mathbb R}}
\def\zed{{\mathbb Z}}
\def\hyp{{\mathbb H}}
\def\hs{{\hat \sigma}}
\def\gmn{{\mathbb G_m^N}}
\def\qbar{{\overline{\mathbb Q}}}
\def\eps{{\varepsilon}}
\def\vek{{\varepsilon_k}}
\def\ahat{{\hat \alpha}}
\def\bhat{{\hat \beta}}
\def\gt{{\tilde \gamma}}
\def\h{{\tfrac12}}
\def\ba{{\boldsymbol a}}
\def\be{{\boldsymbol e}}
\def\bei{{\boldsymbol e_i}}
\def\bc{{\boldsymbol c}}
\def\bm{{\boldsymbol m}}
\def\bk{{\boldsymbol k}}
\def\bi{{\boldsymbol i}}
\def\bl{{\boldsymbol l}}
\def\bq{{\boldsymbol q}}
\def\bu{{\boldsymbol u}}
\def\bt{{\boldsymbol t}}
\def\bs{{\boldsymbol s}}
\def\bv{{\boldsymbol v}}
\def\bw{{\boldsymbol w}}
\def\bx{{\boldsymbol x}}
\def\bX{{\boldsymbol X}}
\def\bz{{\boldsymbol z}}
\def\bwy{{\boldsymbol y}}
\def\bg{{\boldsymbol g}}
\def\bY{{\boldsymbol Y}}
\def\bL{{\boldsymbol L}}
\def\baa{{\boldsymbol\alpha}}
\def\bb{{\boldsymbol\beta}}
\def\bet{{\boldsymbol\eta}}
\def\bxi{{\boldsymbol\xi}}
\def\bo{{\boldsymbol 0}}
\def\bol{{\boldsymbol 1}_L}
\def\ep{\varepsilon}
\def\p{\boldsymbol\varphi}
\def\q{\boldsymbol\psi}
\def\rank{\operatorname{rank}}
\def\aut{\operatorname{Aut}}
\def\lcm{\operatorname{lcm}}
\def\sgn{\operatorname{sgn}}
\def\spn{\operatorname{span}}
\def\md{\operatorname{mod}}
\def\Norm{\operatorname{Norm}}
\def\dim{\operatorname{dim}}
\def\det{\operatorname{det}}
\def\Vol{\operatorname{Vol}}
\def\rk{\operatorname{rk}}

\section{Introduction}
Let 
\begin{equation}
\label{form}
F(\bX,\bY) = \sum_{i=1}^N \sum_{j=1}^N f_{ij} X_i Y_j
\end{equation}
be a symmetric bilinear form in $N \geq 2$ variables with coefficients in a number field $K$. We will also write
$$F(\bX) = F(\bX,\bX)$$
for the associated quadratic form and $F=(f_{ij})_{1 \leq i,j \leq N}$ for the symmetric $N \times N$ coefficient matrix of $F$. We say that the quadratic form $F$ is {\it isotropic} over $K$ if it has a non-trivial zero with coordinates in $K$. A classical theorem of Cassels \cite{cassels:small} states that if $K=\que$ and $F$ is isotropic over $\que$, then there exists $\boldsymbol 0 \neq \bx \in \que^N$ such that $F(\bx)=0$ and
\begin{equation}
\label{cassels}
H(\bx) \ll_N \H(F)^{\frac{N-1}{2}},
\end{equation}
for appropriate notions of heights $H$ and $\H$ to be defined below and an explicit constant. Cassels' theorem has been generalized to any number field $K$ by Raghavan \cite{raghavan}; the bound remained the same as in (\ref{cassels}), except that the explicit constant in the upper bound now depends on the number field as well. 

Further generalizations and extensions of Cassels' theorem have been considered by a number of authors. A few words of notation are required before we can review some of them. Let $Z \subseteq K^N$ be an $L$-dimensional subspace, $1 \leq L \leq N$, then $F$ is defined on $Z$, and we write $(Z,F)$ for the corresponding symmetric bilinear space. A subspace $W$ of $(Z,F)$ is called totally isotropic if for all $\bx,\bwy \in W$, $F(\bx,\bwy)=0$. All maximal totally isotropic subspaces of $(Z,F)$ have the same dimension. It is called the Witt index of $(Z,F)$ and we denote it by $k$.  A subspace $U$ of $(Z,F)$ is called regular if for each $\boldsymbol 0 \neq \bx \in U$ there exists $\bwy \in U$ so that $F(\bx,\bwy) \neq 0$. For each subspace $U$ of $(Z,F)$ we define $U^{\perp} = \{ \bx \in Z : F(\bx, \bwy) = 0\ \forall\ \bwy \in U \}$. If two subspaces $U_1$ and $U_2$ of $(Z,F)$ are orthogonal, we write $U_1 \perp U_2$ for their orthogonal sum. If $U$ is a regular subspace of $(Z,F)$, then $Z = U \perp U^{\perp}$ and $U \cap U^{\perp} = \{\boldsymbol 0\}$.

If $F$ is defined over $\que$, and the bilinear space $(\que^N,F)$ has nonzero Witt index, then Schlickewei \cite{schlickewei} (see also Schmidt and Schlickewei \cite{schmidt:schlickewei}) proved the existence of a maximal totally isotropic subspace of $(\que^N,F)$ of bounded height. This result has been generalized to an arbitrary number field $K$ by Vaaler \cite{vaaler:smallzeros}. In particular, Vaaler proved that if an $L$-dimensional bilinear space $(Z,F)$ over $K$ has Witt index $k \geq 1$, then there exists a maximal totally isotropic subspace $V$ of $(Z,F)$ such that
\begin{equation}
\label{v_bound}
H(V) \ll_{K,L,k} \H(F)^{\frac{L-k}{2}} H(Z).
\end{equation}

The main goal of this paper is to prove a theorem analogous to Vaaler's over $\qbar$, the algebraic closure of $\que$. We use a new method to prove such a result, since Vaaler's argument relies on the fact that the number of subspaces of given dimension and explicitly bounded height over a number field is finite, which is no longer true over $\qbar$. From now on let $(Z,F)$ be an $L$-dimensional regular bilinear space over $\qbar$, $1 \leq L \leq N$. It is a well known fact (see for instance \cite{rylands}) that Witt index of $(Z,F)$ in this case is $\left[ \frac{L}{2} \right]$. We can now state the main result of this paper.

\begin{thm} \label{vaaler} Let $F$ be a quadratic form in $N$ variables as above. Let $Z \subseteq \qbar^N$ be an $L$-dimensional subspace, $1 \leq L \leq N$, so that the quadratic space $(Z,F)$ is regular. Let $k=\left[ \frac{L}{2} \right]$ be the Witt index of $(Z,F)$. There exists a maximal totally isotropic subspace $V$ of $(Z,F)$ with
\begin{equation}
\label{vaaler_bound_even}
H(V) \leq 24 \times 2^{\frac{k-1}{4}} 3^{\frac{k^2(k+1)^2}{4}} \H(F)^{\frac{k^2}{2}} H(Z)^{\frac{k^2+k+2}{2k}},
\end{equation}
if $L$ is even, and
\begin{equation}
\label{vaaler_bound_odd}
H(V) \leq 3^{2k(k+1)^3} \H(F)^{k^2} H(Z)^{\frac{4k}{3}},
\end{equation}
if $L$ is odd.
\end{thm}

Of course given an $L$-dimensional regular bilinear space $(Z,F)$ which is defined over a number field $K$, it is possible to find an extension $E$ of $K$ large enough so that $(Z,F)$ has Witt index $k=\left[ \frac{L}{2} \right]$ over $E$, and then apply Vaaler's theorem with bound (\ref{v_bound}) to it. The constant in (\ref{v_bound}), however, will depend on the discriminant of $E$, which can be quite large. In this case Theorem \ref{vaaler} can produce a stronger bound than (\ref{v_bound}).

Theorem \ref{vaaler} is a statement in the general spirit of ``absolute'' results, in particular it parallels the development of the problem about small-height solutions for a system of homogeneous linear equations, ordinarily known under the name of Siegel's lemma. A version of Siegel's lemma over a number field $K$ asserting the existence of a small-height basis for a subspace of $K^N$ has been proved by Bombieri and Vaaler \cite{vaaler:siegel}. Roy and Thunder \cite{absolute:siegel} proved a version of Siegel's lemma over $\qbar$. More specifically, here is a slightly simplified formulation of the ``absolute'' Siegel's lemma.

\begin{thm}[\cite{absolute:siegel}] \label{roy:thunder} Let $Z \subseteq \qbar^N$ be an $L$-dimensional subspace, $1 \leq L < N$. Then there exists a basis $\bx_1,...,\bx_L$ for $Z$ over $\qbar$ such that
\begin{equation}
\label{Siegel:3}
\prod_{i=1}^L H(\bx_i) \leq \prod_{i=1}^L h(\bx_i) \leq 3^{\frac{L(L-1)}{2}} H(Z).
\end{equation}
\end{thm}

\noindent
In (\ref{Siegel:3}), $h$ stand for inhomogeneous height on vectors to be defined below. Notice that an important common feature distinguishing Theorem \ref{vaaler} and Theorem \ref{roy:thunder} from their number field analogues is that the constants in the upper bounds bear no dependence on any number field.

The proof of Theorem \ref{vaaler} is split into two cases: $L$ is even and $L$ is odd. The argument in the odd case is essentially a reduction to the even case. In the even case we argue by induction on $k=L/2$, the Witt index. We apply the induction hypothesis to the bilinear space $(Z_1,F)$, where $Z_1$ is a codimension two subspace of $Z$ of bounded height guaranteed by the ``absolute'' Siegel's lemma of Roy and Thunder \cite{absolute:siegel}, and hence has Witt index $k-1$. This way we obtain a small-height maximal totally isotropic subspace $U$ of $(Z_1,F)$ guaranteed by the induction hypothesis, and consider its orthogonal dual $W$ in $Z$. We then prove that the intersection of the projective space over $W$ with the quadratic projective variety defined by $F$ over $\qbar$ is a projective intersection cycle whose affine support is a union of two maximal totally isotropic subspaces of $(Z,F)$. Product of their heights can be bounded using a version of arithmetic Bezout's theorem due to Bost, Gillet, and Soul\'e \cite{bgs}.
\smallskip

This paper is structured as follows. In section~2 we set the notation and define the height functions. In section~3 we review a few technical lemmas on properties of heights. In section~4 we prove Theorem \ref{vaaler}. In section~5 we use Theorem \ref{vaaler} to prove an effective version of Witt decomposition theorem for a bilinear space over $\qbar$. We also explain how our results can be extended to bilinear spaces with singular points. In section~6 we derive some related results for quadratic spaces over $\qbar$, including an ``orthogonal'' version of Siegel's lemma and an effective version of Cartan-Dieudonn{\'e} theorem; these are direct analogues of results of \cite{me:witt} over a number field, methods of proof are the same.

\section{Notation and heights}
We start with some notation. Let $K$ be a number field of degree $d$ over $\que$, $O_K$ its ring of integers, and $M(K)$ its set of places. For each place $v \in M(K)$ we write $K_v$ for the completion of $K$ at $v$ and let $d_v = [K_v:\que_v]$ be the local degree of $K$ at $v$, so that for each $u \in M(\que)$
\begin{equation}
\sum_{v \in M(K), v|u} d_v = d.
\end{equation}

\noindent
For each place $v \in M(K)$ we define the absolute value $\|\ \|_v$ to be the unique absolute value on $K_v$ that extends either the usual absolute value on $\real$ or $\cee$ if $v | \infty$, or the usual $p$-adic absolute value on $\que_p$ if $v|p$, where $p$ is a rational prime. We also define the second absolute value $|\ |_v$ for each place $v$ by $|a|_v = \|a\|_v^{d_v/d}$ for all $a \in K$. Then for each non-zero $a \in K$ the {\it product formula} reads
\begin{equation}
\label{product_formula}
\prod_{v \in M(K)} |a|_v = 1.
\end{equation} 

\noindent
We extend absolute values to vectors by defining the local heights. For each $v \in M(K)$ define a local height $H_v$ on $K_v^N$ by
\[ H_v(\bx) = \left\{ \begin{array}{ll}
\max_{1 \leq i \leq N} |x_i|_v & \mbox{if $v \nmid \infty$} \\
\left( \sum_{i=1}^N \|x_i\|_v^2 \right)^{d_v/2d} & \mbox{if $v | \infty$}
\end{array}
\right. \]
for each $\bx \in K_v^N$. We define the following global height function on $K^N$:
\begin{equation}
H(\bx) = \prod_{v \in M(K)} H_v(\bx),
\end{equation}
for each $\bx \in K^N$. Notice that due to the normalizing exponent $1/d$, our global height function is absolute, i.e. for points over $\qbar$ its value does not depend on the field of definition. This means that if $\bx \in \qbar^N$ then $H(\bx)$ can be evaluated over any number field containing the coordinates of $\bx$.

We also define an {\it inhomogeneous} height function on vectors by
 \begin{equation}
h(\bx) = H(1,\bx),
\end{equation}
hence $h(\bx) \geq H(\bx)$ for each $\bx \in \qbar^N$. A basic property of heights that we will use states that for $a_1,...,a_L \in \qbar$ and $\bx_1,...,\bx_L \in \qbar^N$,
\begin{equation}
\label{sum_height}
H \left( \sum_{i=1}^L a_i \bx_i \right) \leq h \left( \sum_{i=1}^L a_i \bx_i \right) \leq H(\ba) \prod_{i=1}^L h(\bx_i),
\end{equation}
where $\ba = (a_1,...,a_L) \in \qbar^L$.
\smallskip

We can extend height to polynomials in the following way: for every $N \geq 1$, if $G(X_1,\dots,X_N) \in K[X_1,\dots,X_N]$ we will write $H_v(G)$ and $H(G)$ for the local and global height of the coefficient vector of $G$, respectively. It is convenient to also introduce a slightly different height $\H$ for our quadratic form $F$: we define $\H_v(F)$ and $\H(F)$ to be the local and global heights, respectively, of the symmetric matrix $(f_{ij})_{1 \leq i,j \leq N}$ viewed as a vector in $K^{N^2}$. It is then easy to see that
\[ H_v(F) \leq \left\{ \begin{array}{ll}
\H_v(F) & \mbox{if $v \nmid \infty$} \\
2^{\frac{d_v}{2d}} \H_v(F) & \mbox{if $v | \infty$},
\end{array}
\right. \]
and so $H(F) \leq \sqrt{2} \H(F)$.
\smallskip

We also define height on matrices, which is the same as height function on subspaces of $\qbar^N$. Let $V \subseteq \qbar^N$ be a subspace of dimension $J$, $1 \leq J \leq N$, defined over a number field $K$. Choose a basis $\bx_1,...,\bx_J$ for $V$ over $K$, and write $X = (\bx_1\ ...\ \bx_J)$ for the corresponding $N \times J$ basis matrix. Then 
$$V = \{ X \bt : \bt \in \qbar^J \}.$$
On the other hand, there exists an $(N-J) \times N$ matrix $A$ with entries in $K$ such that 
$$V = \{ \bx \in \qbar^N : A \bx = 0 \}.$$
Let $\I$ be the collection of all subsets $I$ of $\{1,...,N\}$ of cardinality $J$. For each $I \in \I$ let $I'$ be its complement, i.e. $I' = \{1,...,N\} \setminus I$, and let $\I' = \{ I' : I \in \I\}$. Then 
$$|\I| = \binom{N}{J} = \binom{N}{N-J} = |\I'|.$$
For each $I \in \I$, write $X_I$ for the $J \times J$ sub-matrix of $X$ consisting of all those rows of $X$ which are indexed by $I$, and $_{I'} A$ for the $(N-J) \times (N-J)$ sub-matrix of $A$ consisting of all those columns of $A$ which are indexed by $I'$. By the duality principle of Brill-Gordan \cite{gordan:1} (also see Theorem 1 on p. 294 of \cite{hodge:pedoe}), there exists a non-zero constant $\gamma \in K$ such that
\begin{equation}
\label{duality}
\det (X_I) = (-1)^{\varepsilon(I')} \gamma \det (_{I'} A),
\end{equation}
where $\varepsilon(I') = \sum_{i \in I'} i$. Define the vectors of {\it Grassmann coordinates} of $X$ and $A$ respectively to be 
$$Gr(X) = (\det (X_I))_{I \in \I} \in K^{|I|},\ \ Gr(A) = (\det (_{I'} A))_{I' \in \I'} \in K^{|I'|}.$$
Define 
$$H(X) =  H(Gr(X)),\ \ H(A) = H(Gr(A)),$$
and so by (\ref{duality}) and (\ref{product_formula})
$$H(X) = H(A).$$
Define height of $V$ denoted by $H(V)$ to be this common value. Hence the height of a matrix is the height of its row (or column) space, which is equal to the height of its null-space. Also notice that $Gr(X)$ can be identified with $\bx_1 \wedge\ ...\ \wedge \bx_J$, where $\wedge$ stands for the wedge product, viewed under the canonical lexicographic embedding into $K^{\binom{N}{J}}$. Therefore we can also write
$$H(V) = H(\bx_1 \wedge\ ...\ \wedge \bx_J).$$
This definition is legitimate, since it does not depend on the choice of the basis for $V$: let $\bwy_1,...,\bwy_J$ be another basis for $V$ over $K$, then there exists $C \in GL_N(K)$ such that $\bwy_i = C \bx_i$ for each $1 \leq i \leq J$, and so
\begin{eqnarray*}
H(\bwy_1 \wedge\ ...\ \wedge \bwy_J) & = & H(C \bx_1 \wedge\ ...\ \wedge C \bx_J) \\
& = & \left( \prod_{v \in M(K)} |\det(C)|_v \right) H(\bx_1 \wedge\ ...\ \wedge \bx_J) \\
& = & H(\bx_1 \wedge\ ...\ \wedge \bx_J),
\end{eqnarray*}

Finally, for a point $\bz = (z_1,...,z_N) \in \qbar^N$, we write $\deg_K(\bz)$ to mean the degree of the extension $K(z_1,...,z_N)$ over $K$, i.e. $\deg_K(\bz) = [K(z_1,...,z_N):K]$.

\section{Preliminary lemmas}

Here we present some technical lemmas that we will use. The first one is a consequence of Laplace's expansion, and can be found as Lemma 4.7 of \cite{absolute:siegel} (also see pp. 15-16 of \cite{vaaler:siegel}).

\begin{lem} \label{wedge} Let $X$ be a $N \times J$ matrix over $\qbar$ with column vectors $\bx_1,...,\bx_J$. Then
\begin{equation}
\label{prod_1}
H(X) = H(\bx_1 \wedge \bx_1\ ...\ \wedge \bx_J) \leq \prod_{i=1}^J H(\bx_i).
\end{equation}
More generally, if the $N \times J$ matrix $X$ can be partitioned into blocks as $X = (X_1\ X_2)$, then
 \begin{equation}
\label{prod_2}
H(X) \leq H(X_1) H(X_2).
\end{equation}
\end{lem}

The next one is an obvious adaptation of Lemma 2.3 of \cite{me:witt} over~$\qbar$. 

\begin{lem} \label{mult} Let $X$ be a $N \times J$ matrix over $\qbar$ with column vectors $\bx_1,...,\bx_J$, and let $F$ be a symmetric bilinear form in $N$ variables, as above (we also write $F$ for its $N \times N$ coefficient matrix). Then
\begin{equation}
\label{prod_3}
H(F X) \leq \H(F)^J \prod_{i=1}^J H(\bx_i).
\end{equation}
\end{lem}

The following well known fact is an immediate corollary of Theorem 1 of \cite{vaaler:struppeck} adapted over $\qbar$.

\begin{lem} \label{intersection} Let $U_1$ and $U_2$ be subspaces of $\qbar^N$. Then
$$H(U_1 \cap U_2) \leq H(U_1) H(U_2).$$
\end{lem}

The following simple lemma can be viewed as an analogue of Cassels' bound (\ref{cassels}) over $\qbar$. It is a special case of Proposition 3.1 of \cite{me:bezout}; we include a proof here for the purposes of self-containment.

\begin{lem} \label{qz} Let $F$ be a quadratic form in $N$ variables as above. Then there exists $\boldsymbol 0 \neq \bx \in \qbar^N$ such that $F(\bx) = 0$, and 
\begin{equation}
\label{qz:bound}
H(\bx) \leq 2 \sqrt{\H(F)}.
\end{equation}
\end{lem}

\proof
If $F$ is identically zero, then we are done. So assume $F$ is non-zero. Write $\be_1,...,\be_N$ for the standard basis vectors for $\qbar^N$ over $\qbar$. Assume that for some $1 \leq i \leq N$, $\deg_{X_i} F < 2$, then it is easy to see that $F(\be_i)=0$, and $H(\be_i) = 1$. If $N>2$, let
$$F_1(X_1,X_2) = F(X_1,X_2,0,...,0),$$
and a point $\bx = (x_1,x_2) \in \qbar^2$ is a zero of $F_1$ if and only if $(x_1,x_2,0,...,0)$ is a zero of $F$, and $H(x_1,x_2) = H(x_1,x_2,0,...,0)$. In particular, if $F_1(X_1,X_2) = 0$, then $F(\be_1) = 0$. Hence we can assume that $N=2$, $F(X_1,X_2) \neq 0$, and $\deg_{X_1} F = \deg_{X_2} F = 2$. Write
$$F(X_1,X_2) = f_{11} X_1^2 + 2 f_{12} X_1 X_2 + f_{22} X_2^2,$$
where $f_{11}, f_{22} \neq 0$. Let
$$g(X_1) = F(X_1,1) = f_{11} X_1^2 + 2f_{12} X_1 + f_{22},$$
be a quadratic polynomial in one variable. Notice that since coefficients of $g$ are those of $F$, we have $H(g) = H(F) \leq \sqrt{2} \H(F)$. By Lemma 2 of \cite{vaaler:pinner} (see also Lemma 2.1 of \cite{me:bezout}), there must exist $\alpha \in \qbar$ such that $g(\alpha) = 0$, and
$$H(\alpha,1) \leq \sqrt{2 H(g)} \leq 2 \sqrt{\H(F)}.$$
Taking $\bx = (\alpha,1)$ completes the proof.
\endproof

Next we present a lemma on effective decomposition of a quadratic space into regular and singular components. This is an adaptation of Lemma 3.2 of \cite{me:witt} over $\qbar$, although the inequality (\ref{regular2}) is slightly weaker than its number field analogue; this, however, makes essentially no difference for our purposes.

\begin{lem} \label{reg:decompose} Let $F$ have rank $r$ on $Z$, and assume that $1 \leq r < L$. Then the bilinear space $(Z,F)$ can be represented as
\begin{equation}
\label{regular1}
Z = Z^{\perp} \perp W,
\end{equation}
where $Z^{\perp} = \{ \bz \in Z : F(\bz,\bx) = 0\ \forall\ \bx \in Z \}$ is the $(L-r)$-dimensional singular component, and $W$ is a regular subspace of $Z$, with $\dim_{\qbar} W = r$ and
\begin{equation}
\label{regular2}
H(Z^{\perp}) \leq 3^{\frac{L(L-1)}{2}} \H(F)^r H(Z)^2,
\end{equation}
and
\begin{equation}
\label{regular3}
H(W) \leq 3^{\frac{L(L-1)}{2}} H(Z).
\end{equation}
\end{lem}

\proof
Let $\bx_1, \dots ,\bx_L$ be the basis for $Z$ guaranteed by Theorem \ref{roy:thunder}, and write $X = (\bx_1 \dots \bx_L)$ for the corresponding $N \times L$ basis matrix. Notice that
$$Z^{\perp} = \Nn(FX) \cap Z,$$
where $\Nn(FX) = \{ \bz \in \qbar^N : \bz FX = 0 \}$ is the null-space of the matrix $FX$. Since the matrix $FX = (F\bx_1 \dots F\bx_L)$ has rank $r < L$, only $r$ of its columns can be linearly independent. In other words, there exist $\bx_{i_1},...,\bx_{i_r} \in \{ \bx_1,...,\bx_L \}$ such that the $N \times r$ matrix $FX' = (F\bx_{i_1} \dots F\bx_{i_L})$ has rank $r$, and $\Nn(FX) = \{ \bz \in \qbar^N : \bz FX' = 0 \}$. Then, combining Lemmas \ref{mult} and \ref{intersection} with Theorem \ref{roy:thunder}, we obtain
\begin{eqnarray*}
H(Z^{\perp}) & \leq & H(\Nn(FX)) H(Z) = H(FX') H(Z) \leq H(Z) \H(F)^r \prod_{j=1}^r  H(\bx_{i_j}) \\
& \leq & H(Z) \H(F)^r \prod_{i=1}^L H(\bx_i) \leq 3^{\frac{L(L-1)}{2}} \H(F)^r H(Z)^2,
\end{eqnarray*}
which is precisely (\ref{regular2}). The proof of (\ref{regular3}) is identical to the argument in the proof of Lemma 3.2 of \cite{me:witt}, but we use Theorem \ref{roy:thunder} instead of Bombieri-Vaaler version of Siegel's lemma.
\endproof

Finally, we will need the following result on the existence of a small-height vector in a subspace $Z \subseteq \qbar^N$ at which the quadratic form $F$ does not vanish, satisfying one additional condition. An {\it isometry} of the quadratic space $(Z,F)$ is an isomorphism $\sigma : Z \rightarrow Z$ such that $F(\sigma(\bx), \sigma(\bwy)) = F(\bx, \bwy)$ for all $\bx, \bwy \in Z$. It is easy to see that isometries of $(Z,F)$ form a group under function composition, which we denote by $\OO(Z,F)$; isometries will be discussed in more details in section~6. The following is a direct analogue of Lemma 5.2 of \cite{me:witt} over $\qbar$.

\begin{lem} \label{nonzero} Let $(Z,F)$ be an $L$-dimensional regular bilinear space in $N$ variables over $\qbar$, as above, and let $\sigma \in \OO(Z,F)$. There exists an anisotropic vector $\bwy$ in $Z$ such that $\sigma(\bwy) \pm \bwy$ is also anisotropic for some choice of $\pm$, and
\begin{equation}
\label{anis}
H(\bwy) \leq h(\bwy) \leq 2 \sqrt{L}\ 3^{\frac{(L+2)(L-1)}{4}} H(Z)^{\frac{L+2}{2L}}.
\end{equation}
\end{lem}

\proof
The argument is identical to that in the proof of Lemma 5.2 of \cite{me:witt} with Bombieri-Vaaler version of Siegel's lemma \cite{vaaler:siegel} replaced with the absolute version of Roy and Thunder (Theorem \ref{roy:thunder}).
\endproof

We are now ready to proceed.

\section{Proof of Theorem \ref{vaaler}}
Let $F$ be a quadratic form in $N \geq 2$ variables, as above. Let $Z \subseteq \qbar^N$ be an $L$-dimensional subspace, $1 \leq L \leq N$, such that the bilinear space $(Z,F)$ is regular. We start by proving a lemma about the existence of a small-height zero of $F$ in $Z$.

\begin{lem} \label{zero} Let $2 \leq L \leq N$. There exists $\boldsymbol 0 \neq \bwy \in Z$ such that $F(\bwy)=0$ and
\begin{equation}
\label{zero_eq}
H(\bwy) \leq 8 \times 3^{2(L-1)} \H(F)^{\frac{1}{2}} H(Z)^{\frac{4}{L}}.
\end{equation}
\end{lem}

\proof
 Let $\bz_1,\dots, \bz_L$ be the basis for $Z$ guaranteed by Theorem \ref{roy:thunder} ordered so that
$$H(\bz_1) \leq \dots \leq H(\bz_L).$$
Then, by (\ref{Siegel:3}),
\begin{equation}
\label{z1}
H(\bz_1) H(\bz_2) \leq \left( 3^{\frac{L(L-1)}{2}} H(Z) \right)^{\frac{2}{L}} = 3^{L-1} H(Z)^{\frac{2}{L}}.
\end{equation}
We will now construct $a_1,a_2 \in \qbar$ such that $\bwy = a_1 \bz_1 + a_2 \bz_2$ is a zero of $F$. In other words, we want
\begin{equation}
\label{newform}
0 = F(\bwy) = F(\bz_1) a_1^2 + 2F(\bz_1,\bz_2) a_1 a_2 + F(\bz_2) a_2^2 = G(a_1,a_2).
\end{equation}
The right hand side of (\ref{newform}) is a quadratic form $G$ in the variables $a_1,a_2$ with coefficients $F(\bz_1), 2F(\bz_1,\bz_2),F(\bz_2)$. By Lemma \ref{qz}, there must exist such a pair $(a_1,a_2)$ with
\begin{equation}
\label{q0}
H(a_1,a_2) \leq 2 \sqrt{\H(G)}.
\end{equation}
Let $E$ be the field extension generated over $K$ by coefficients of $G$. By (2.6) of \cite{vaaler:smallzeros2}, for each $v \in M(E)$ and each $\baa_1,\baa_2 \in E_v^N$
$$|F(\baa_1,\baa_2)|_v \leq \H_v(F) H_v(\baa_1) H_v(\baa_2).$$
Therefore, if $v \nmid \infty$, we have
\begin{eqnarray}
\label{q1}
\H_v(G) & \leq & \max \{ |F(\bz_1)|_v, |2|_v |F(\bz_1,\bz_2)|_v, |F(\bz_2)|_v \} \nonumber \\
& \leq & \H_v(F) \max \{ H_v(\bz_1)^2, H_v(\bz_1)H_v(\bz_2), H_v(\bz_2)^2 \} \\
& \leq & \H_v(F) \max \{ 1, H_v(\bz_1) \}^2 \max \{ 1, H_v(\bz_2) \}^2. \nonumber
\end{eqnarray}
If $v | \infty$, then
\begin{eqnarray}
\label{q2}
\H_v(G)^{\frac{2d}{d_v}} & \leq & \|F(\bz_1)\|^2_v + 2 \|F(\bz_1,\bz_2)\|^2_v + \|F(\bz_2)\|^2_v \leq \H_v(F)^{\frac{2d}{d_v}} \times \nonumber \\
& \times & \left( H_v(\bz_1)^{\frac{4d}{d_v}} + 2 \left( H_v(\bz_1) H_v(\bz_2) \right)^{\frac{2d}{d_v}} + H_v(\bz_2)^{\frac{4d}{d_v}} \right) \\
& \leq & \H_v(F)^{\frac{2d}{d_v}} \left( 1 + H_v(\bz_1)^{\frac{2d}{d_v}} \right)^2 \left( 1 + H_v(\bz_2)^{\frac{2d}{d_v}} \right)^2. \nonumber
\end{eqnarray}
Combining (\ref{q0}) with (\ref{q1}) and (\ref{q2}), we see that
\begin{equation}
\label{q3}
H(a_1,a_2) \leq 2 \H(F)^{\frac{1}{2}} h(\bz_1) h(\bz_2).
\end{equation}
Combining (\ref{q3}) and (\ref{sum_height}), we see that there exists a zero of $F$ of the form $\bwy = a_1 \bz_1 + a_2 \bz_2 \in Z$ so that 
\begin{equation}
\label{q4}
H(\bwy) \leq 2 \H(F)^{\frac{1}{2}} h(\bz_1)^2 h(\bz_2)^2.
\end{equation}
Notice that for each $l=1,2$, $\bz_l$ is a non-zero vector. If $\bz_l$ has just one non-zero coordinate, let it for instance be $i$-th coordinate, then clearly we can take $\bz_l$ to be $i$-th standard basis vector $\bei$, and so $h(\bz_l) = \sqrt{2}$. If $\bz_l$ has more than one non-zero coordinates, then we can assume without loss of generality that at least one of them is equal to $1$, and then it is easy to see that $h(\bz_l) \leq \sqrt{2} H(\bz_l)$. Therefore (\ref{q4}) can be rewritten as
\begin{equation}
\label{q5}
H(\bwy) \leq 8 \H(F)^{\frac{1}{2}} \left( H(\bz_1) H(\bz_2) \right)^2 \leq 8 \times 3^{2(L-1)} \H(F)^{\frac{1}{2}} H(Z)^{\frac{4}{L}},
\end{equation}
where the last inequality follows by (\ref{z1}). This completes the proof.
\endproof

We are now ready to prove the theorem.
\bigskip

{\it Proof of Theorem \ref{vaaler}.}
First suppose that $L$ is even, say $L=2k$ for some integer $k \geq 1$. We argue by induction on $k$, the Witt index of $(Z,F)$. Suppose that $k=1$. Let $\bwy \in Z$ be the point guaranteed by Lemma \ref{zero}. Then
\begin{equation}
\label{q6}
H(\bwy) \leq 72\ \H(F)^{\frac{1}{2}} H(Z)^{2},
\end{equation}
and $V = \spn_{\qbar} \{\bwy\}$ is the desired maximal totally isotropic subspace of $(Z,F)$.

Next suppose $k>1$. Let $\bz_1,\dots, \bz_L$ be the basis for $Z$ guaranteed by Theorem \ref{roy:thunder} ordered so that
$$H(\bz_1) \leq \dots \leq H(\bz_L).$$
Let $Z_1 = \spn_{\qbar} \{\bz_1,...,\bz_{L-2}\}$, so that
\begin{equation}
\label{ind_sub}
H(Z_1) \leq \prod_{i=1}^{L-2} H(\bz_i) \leq \left( 3^{\frac{L(L-1)}{2}} H(Z) \right)^{\frac{L-2}{L}} \leq 3^{(2k-1)(k-1)} H(Z)^{\frac{k-1}{k}}.
\end{equation}
Let $U$ be the maximal totally isotropic subspace of bounded height of $(Z_1,F)$ guaranteed by the induction hypothesis, so that 
\begin{eqnarray}
\label{ind_hyp}
H(U) & \leq & 24 \times 2^{\frac{k-2}{4}} 3^{\frac{k^2(k-1)^2}{4}} \H(F)^{\frac{(k-1)^2}{2}} H(Z_1)^{\frac{k^2-k+2}{2(k-1)}} \nonumber \\
& \leq & 24 \times 2^{\frac{k-2}{4}} 3^{\frac{(2k-1)(k^2-k+2)}{2}+\frac{k^2(k-1)^2}{4}} \H(F)^{\frac{(k-1)^2}{2}} H(Z)^{\frac{k^2-k+2}{2k}},
\end{eqnarray}
by (\ref{ind_sub}). Since $\dim_{\qbar} (U) = k-1$, let $\bx_1,...,\bx_{k-1}$ be the basis for $U$ guaranteed by Theorem \ref{roy:thunder}, and let $X = (\bx_1\ ...\ \bx_{k-1})$ be the corresponding $N \times (k-1)$ basis matrix. Define a $(k+1)$-dimensional subspace of $Z$
$$W = \{ \bwy \in Z : \bwy FX = 0 \} = Z \cap \{ \bwy \in \qbar : \bwy FX = 0 \}.$$
Combining Lemma \ref{mult}, Lemma \ref{intersection}, and Theorem \ref{roy:thunder}, we see that
\begin{eqnarray}
\label{H_W}
H(W) \leq H(Z) H(F X) & \leq & H(Z) \H(F)^{k-1} \prod_{i=1}^{k-1} H(\bx_i) \nonumber \\
& \leq & \left\{ 3^{\frac{k}{2}} \H(F) \right\}^{k-1} H(U) H(Z).
\end{eqnarray}
It is easy to see that $U \subset W$. Let $\bw_1,...,\bw_{k+1}$ be a basis for $W$. Then at least two of these basis vectors are not in $U$, we can assume without loss of generality that these are $\bw_1$ and $\bw_2$. Since $W=\spn_{\qbar} \{U,\bw_1,\bw_2\}$, and $\dim_{\qbar}(W)=\dim_{\qbar}(U)+2$, it must be true that $U \cap \spn_{\qbar} \{\bw_1,\bw_2\} = \{\boldsymbol 0\}$. Consider a binary quadratic form $G$ in two variables $a,b$ given by
$$G(a,b) = F(\bw_1) a^2 + 2F(\bw_1,\bw_2) a b + F(\bw_2) b^2.$$
Let $a=1$, then there exist $b_1,b_2 \in \qbar$ such that $G(1,b_1)=G(1,b_2)=0$, i.e. the zero-set of $G$ consists, up to multiplicity, of two projective points. Let $\bwy_i = \bw_1 + b_i \bw_2$ for each $i=1,2$, and so $F(\bwy_1)=F(\bwy_2)=0$ and $\bwy_1,\bwy_2 \notin U$. Define $V_i=\spn_{\qbar} \{U,\bwy_i\}$ for each $i=1,2$, then $V_1,V_2 \subset W$ are maximal totally isotropic subspaces of $(Z,F)$. Now suppose that $\bx \in W$, then $\bx = \bx' + \beta_1 \bw_1 + \beta_2 \bw_2$, where $\bx' \in U$. Therefore 
$$F(\bx) = F(\beta_1 \bw_1 + \beta_2 \bw_2) = G(\beta_1,\beta_2).$$
Hence $F(\bx) = 0$ if and only if $(\beta_1,\beta_2)$ is a multiple of either $(1,b_1)$ or $(1,b_2)$. In other words, $\bx \in W$ is such that $F(\bx)=0$ if and only if $\bx \in V_1 \cup V_2$.

Let $\Z_F $ be the projective closure of the affine set $\{ \bx \in \qbar^N : F(\bx) = 0 \}$, i.e. the projective variety defined by $F$ over $\qbar$. Write $\pee(W)$ for the projective space over $W$. Let $\Z_F \centerdot \pee(W)$ be the intersection cycle of these two projective varieties, so that up to multiplicity
$$\Z_F \centerdot \pee(W) = \pee(V_1) + \pee(V_2),$$
i.e. its support is $\pee(V_1) \cup \pee(V_2)$. Notice that from our construction above it is possible that $V_1=V_2$, then we write $\Z_F \centerdot \pee(W) = 2 \pee(V_1)$, so that 2 is the multiplicity of the component $V_1$ in this intersection cycle. In any case, the height of this intersection cycle is defined by
\begin{equation}
\label{hint1}
H(\Z_F \centerdot \pee(W)) = H(V_1) H(V_2).
\end{equation}
See \cite{bgs} for the details on arithmetic intersection theory and heights, keeping in mind that in case of a linear space or a hypersurface their height reduces to an additive height given by $\log (H^d)$ (see (3.1.6) on p. 947 and remark after Proposition 4.1.2 on p. 965 of \cite{bgs}). Then applying a version of Arithmetic Bezout's Theorem presented by Theorem 5.4.4 (i) of \cite{bgs}, we obtain
\begin{equation}
\label{hint2}
H(\Z_F \centerdot \pee(W)) \leq \sqrt{2} H(W)^{2} \H(F).
\end{equation}
Therefore, combining (\ref{hint1}), (\ref{hint2}), and (\ref{H_W}), we have
\begin{eqnarray}
\label{hint3}
\min \{H(V_1),H(V_2)\} & \leq & \sqrt{H(V_1)H(V_2)} = \sqrt{H(\Z_F \centerdot \pee(W))} \nonumber \\
& \leq & 2^{\frac{1}{4}} H(W) \H(F)^{\frac{1}{2}} \leq 2^{\frac{1}{4}} 3^{\frac{k(k-1)}{2}} \H(F)^{k-\frac{1}{2}} H(U) H(Z).
\end{eqnarray}
Write $V$ for $V_i$ with $H(V_i)$ being the smaller of the two, $i=1,2$. Combining (\ref{ind_hyp}) and (\ref{hint3}) yields (\ref{vaaler_bound_even}).
\smallskip

Next suppose that $L$ is odd, say $L=2k+1$ for some $k \geq 0$. We again argue by induction on $k$, which is the Witt index of $(Z,F)$. If $k=0$, then $L=1$, and so $Z=\qbar \bwy$ for some anisotropic vector $\bwy \in \qbar^N$ since $(Z,F)$ is regular. Hence $\{ \boldsymbol 0 \}$ is the maximal totally isotropic subspace of $(Z,F)$, so there is nothing to prove. Assume $k \geq 1$. Let $\bz_1,\dots, \bz_L$ be the basis for $Z$ guaranteed by Theorem \ref{roy:thunder} ordered so that
$$H(\bz_1) \leq \dots \leq H(\bz_L),$$
and define $Z_1 = \spn_{\qbar} \{ \bz_1,...,\bz_{L-1} \}$. By Theorem \ref{roy:thunder}
\begin{equation}
\label{Z1}
H(Z_1) \leq 3^{2k^2} H(Z)^{\frac{2k}{2k+1}}.
\end{equation}
Then $(Z_1,F)$ is a bilinear space over $\qbar$ of dimension $L-1=2k$. Notice that maximal totally isotropic subspaces of $(Z_1,F)$ are also maximal totally isotropic subspaces of $(Z,F)$, hence have dimension $k$. Let
$$Z_1^{\perp} = \{ \bx \in Z_1 : F(\bx,\bz) = 0\ \forall\ \bz \in Z_1 \}$$
be the singular component of $Z_1$, and suppose it has dimension $l$. Then $Z_1 = Z_1^{\perp} \perp W$ for some $(2k-l)$-dimensional subspace $W$ of $Z_1$ such that $(W,F)$ is regular. Witt index of $(W,F)$ is therefore equal to $\left[ \frac{2k-l}{2} \right]$, so let $V$ be a maximal totally isotropic subspace of $W$. Then $Z_1^{\perp} \perp V$ is a maximal totally isotropic subspace of $Z_1$, and so has dimension $k$. On the other hand, $Z_1^{\perp} \cap V = \{ \boldsymbol 0 \}$, hence $\dim_{\qbar} (Z_1^{\perp} \perp V) = l + \left[ \frac{2k-l}{2} \right]$. Therefore we must have
$$l + \left[ \frac{2k-l}{2} \right] = k,$$
which means that either $l=0$ or $l=1$. If $l=0$, then $(Z_1,F)$ is regular, and so by the argument in the even case above there exists a maximal totally isotropic subspace $V$ of $(Z_1,F)$ of bounded height. Moreover, by combining (\ref{vaaler_bound_even}) and (\ref{Z1}), we obtain
\begin{equation}
\label{o1}
H(V) \leq 24 \times 2^{\frac{k-1}{4}} 3^{\frac{k(k+1)^2(k+4)}{4}} \H(F)^{\frac{k^2}{2}} H(Z)^{\frac{k^2+k+2}{2k+1}},
\end{equation}
which is smaller than the bound in (\ref{vaaler_bound_odd}).

Assume $l=1$. Then $(W,F)$ is a regular $(2k-1)$-dimensional bilinear space of Witt index $k-1$ with
\begin{equation}
\label{W}
H(W) \leq 3^{k(4k-1)} H(Z)^{\frac{2k}{2k+1}},
\end{equation}
where this bound is obtained by combining (\ref{regular3}) with (\ref{Z1}). By induction hypothesis, there exists a maximal totally isotropic subspace $U$ of $(W,F)$ with
\begin{equation}
\label{ind_odd}
H(U) \leq 3^{2(k-1)k^3} \H(F)^{(k-1)^2} H(W)^{\frac{4(k-1)}{3}}.
\end{equation}
Notice that $F$ has rank $2k-1$ on $Z_1$. Combining Lemma \ref{reg:decompose} and (\ref{ind_odd}), and applying (\ref{prod_2}) of Lemma \ref{wedge}, we see that the maximal totally isotropic subspace $V = Z_1^{\perp} \perp U$ of $(Z_1,F)$ satisfies
\begin{equation}
\label{odd1}
H(V) \leq H(Z_1^{\perp}) H(U) \leq 3^{k(6k-1)} \H(F)^{2k-1} H(Z)^{\frac{4k}{2k+1}} H(U),
\end{equation}
where the last inequality follows by combining (\ref{regular2}) and (\ref{Z1}). Combining (\ref{W}), (\ref{ind_odd}), and (\ref{odd1}) produces (\ref{vaaler_bound_odd}), and so finishes the proof.
\boxed{ }

\section{Corollaries}

In this section we use Theorem \ref{vaaler} to prove an effective version of Witt decomposition theorem over $\qbar$. First let us recall that a {\it hyperbolic plane} in $(Z,F)$ is a two-dimensional subspace of the form
$$\hyp = \spn_{\qbar} \{ \bx, \bwy \in Z : F(\bx) = F(\bwy) = 0,\ F(\bx,\bwy) = 1 \}.$$
A classical theorem of Witt (1937) for a regular $L$-dimensional bilinear space $(Z,F)$ over $\qbar$ states that there exists a decomposition of $Z$ into an orthogonal direct sum with respect to $F$ of the form
\begin{equation}
\label{witt}
Z = \hyp_1 \perp \dots \perp \hyp_k \perp W,
\end{equation}
where $k=\left[ \frac{L}{2} \right]$ is Witt index of $(Z,F)$, $\hyp_1,\dots,\hyp_k$ are hyperbolic planes, and $W$ is zero if $L=2k$, and is a one-dimensional anisotropic component if $L=2k+1$, i.e. is of the form $\spn_{\qbar} \{\bwy\}$ for some $\bwy \in Z$ such that $F(\bwy) \neq 0$ (see for instance Corollary 5.11 on p.17 of \cite{scharlau}). An effective version of Witt decomposition theorem for bilinear spaces over a number field is proved in \cite{me:witt}. Here we obtain the following effective analogue of Witt's theorem over $\qbar$.

\begin{thm} \label{witt_qbar} Let $(Z,F)$ be a regular $L$-dimensional bilinear space in $N$-variables over $\qbar$. There exists an orthogonal decomposition of $(Z,F)$ as in (\ref{witt}) such that for each $1 \leq i \leq k=\left[ \frac{L}{2} \right]$
\begin{equation}
\label{witt1}
H(\hyp_i) \leq 3^{12 k^4 (k+1) \left(\frac{3}{2}\right)^k} \left\{ \sqrt{k}\ \H(F)^{k^2+1} H(Z)^{\frac{6k+5}{4k+2}} \right\}^{\frac{(k+1)(k+2)}{2} \left(\frac{3}{2}\right)^k},
\end{equation}
and $W=\{\boldsymbol 0\}$ if $L=2k$, or $W=\qbar \bwy$ with
\begin{equation}
\label{witt2}
H(W) = H(\bwy) \leq 2 \sqrt{2k+1}\ 3^{\frac{(2k+3)k}{2}} H(Z)^{\frac{2k+3}{4k+2}},
\end{equation}
if $L=2k+1$.
\end{thm}

\proof
First assume $L=2k$. In this case we prove that there exists a decomposition of $(Z,F)$ of the form
\begin{equation}
\label{dh}
Z = \hyp_1 \perp \dots \perp \hyp_k,
\end{equation}
with
\begin{equation}
\label{even_dec}
H(\hyp_i) \leq 3^{12 k^5 \left(\frac{3}{2}\right)^k} \left\{ \H(F)^{k^2} H(Z) \right\}^{\frac{(k+1)(k+2)}{2} \left(\frac{3}{2}\right)^k},
\end{equation}
for each $1 \leq i \leq k$, which is smaller than the bound in (\ref{witt1}). The argument is identical to the proof of Theorem 4.1 of \cite{me:witt}, except that instead of Bombieri-Vaaler version of Siegel's lemma over a fixed number field (quoted as Theorem 2.1 in \cite{me:witt}) we use the absolute version of our Theorem \ref{roy:thunder} due to Roy and Thunder, and instead of Vaaler's theorem on the existence of a maximal totally isotropic subspace of bounded height in a bilinear space over a fixed number field (quoted as Theorem 3.1 in \cite{me:witt}) we use the absolute version, i.e. the case $L=2k$ of our Theorem \ref{vaaler} given by the bound in (\ref{vaaler_bound_even}). Using the same construction as in the proof of Theorem 4.1 of \cite{me:witt}, we first obtain a hyperbolic plane $\hyp_1 \subseteq Z$ such that
\begin{equation}
\label{hyp1}
H(\hyp_1) \leq 3^{(k+1)^3} \H(F)^{\frac{k}{2}} H(Z)^{\frac{(k+1)(k+2)}{2k^2}},
\end{equation}
which is smaller than the bound in (\ref{even_dec}), and then define 
$$Z_1 = \hyp_1^{\perp} = \{ \bz \in \qbar^N : F(\bz,\bx) = 0\ \forall\ \bx \in \hyp_1 \} \cap Z,$$
so that $\dim_{\qbar} (Z_1) = L-2 = 2(k-1)$, $(Z_1,F)$ is a regular bilinear space of Witt index $k-1$, and $Z = \hyp_1 \perp Z_1$. Combining Lemmas \ref{mult} and \ref{intersection} with (\ref{hyp1}), we obtain
\begin{equation}
\label{zz1}
H(Z_1) \leq H(\hyp_1) H(Z) \H(F)^2 \leq 3^{(k+1)^3} \H(F)^{\frac{k+4}{2}} H(Z)^{\frac{3k^2+3k+2}{2k^2}}.
\end{equation}
We proceed by induction on $k$. If $k=1$, we are done. If $k \geq 2$, by induction hypothesis there must exist a decomposition for $(Z_1,F)$ of the form
$$Z_1 = \hyp_2 \perp \dots \perp \hyp_k,$$
with
\begin{equation}
\label{dec_hyp}
H(\hyp_i) \leq 3^{12 (k-1)^5 \left(\frac{3}{2}\right)^{k-1}} \left\{ \H(F)^{(k-1)^2} H(Z_1) \right\}^{\frac{k(k+1)}{2} \left(\frac{3}{2}\right)^{k-1}},
\end{equation}
for each $2 \leq i \leq k$. Then 
$$Z = \hyp_1 \perp Z_1 = \hyp_1 \perp \hyp_2 \perp \dots \perp \hyp_k,$$
and combining (\ref{zz1}) and (\ref{dec_hyp}) yields (\ref{even_dec}).
\smallskip

Next suppose $L=2k+1$. Then let $\bwy \in Z$ be the anisotropic vector guaranteed by Lemma \ref{nonzero}, and define $W=\spn_{\qbar} \{\bwy\}$, so that by (\ref{anis})
\begin{equation}
\label{anis_comp}
H(W) = H(\bwy) \leq 2 \sqrt{2k+1}\ 3^{\frac{(2k+3)k}{2}} H(Z)^{\frac{2k+3}{4k+2}}.
\end{equation}
Define 
$$Z_1 = W^{\perp} = \{ \bz \in \qbar^N : F(\bz,\bwy) = 0 \} \cap Z,$$
then $\dim_{\qbar} (Z_1) = L-1 = 2k$, and $Z=W \perp Z_1$. Combining Lemmas \ref{mult} and \ref{intersection} with (\ref{anis_comp}), we obtain
\begin{equation}
\label{zz2}
H(Z_1) \leq H(W) H(Z) \H(F) \leq 2 \sqrt{2k+1}\ 3^{\frac{(2k+3)k}{2}} \H(F) H(Z)^{\frac{6k+5}{4k+2}}.
\end{equation}
Suppose there exists $\bx \in Z_1$ such that $F(\bx,\bz)=0$ for all $\bz \in Z_1$. By construction $F(\bx,\bwy)=0$, and hence $F(\bx,\bz)=0$ for all $\bz \in Z$. This contradicts the original assumption that $(Z,F)$ is regular, hence $(Z_1,F)$ must also be regular, and thus its Witt index is equal to $k$. Therefore, by the argument in the even case above, there exists a decomposition into an orthogonal sum of $k$ hyperbolic planes of bounded height like (\ref{dh}) for $Z_1$. Hence we obtained a decomposition as in (\ref{witt}) for $Z$. Combining (\ref{even_dec}) with (\ref{zz2}) yields (\ref{witt1}), and (\ref{witt2}) is precisely (\ref{anis_comp}). 
\endproof

{\it Remark.} Notice that Theorems \ref{vaaler} and \ref{witt_qbar} can both be extended to the case when our bilinear space $(Z,F)$ contains singular points. In this case $Z$ can be written as $Z = Z^{\perp} \perp Z_1$, where $Z^{\perp}$ is the singular component and $(Z_1,F)$ is a regular bilinear space. If $V$ is the maximal totally isotropic subspace of $(Z_1,F)$ of bounded height as guaranteed by Theorem \ref{vaaler}, then $Z^{\perp} \perp V$ is a  maximal totally isotropic subspace of $(Z,F)$, and by Lemma \ref{wedge}
$$H(Z^{\perp} \perp V) \leq H(Z^{\perp}) H(V),$$
where $H(Z^{\perp})$ is bounded by (\ref{regular2}) of Lemma \ref{reg:decompose}. Witt decomposition for $(Z,F)$ in this case will be of the form
$$Z = Z^{\perp} \perp \hyp_1 \perp \dots \perp \hyp_k \perp W,$$
where $k$ is the Witt index of $(Z_1,F)$, and $Z_1 = \hyp_1 \perp \dots \perp \hyp_k \perp W$ is the small-height decomposition for $(Z_1,F)$ guaranteed by Theorem \ref{witt_qbar}.

\section{Related Results}

In this section we derive some additional structural theorems for bilinear spaces over $\qbar$, analogous to those in \cite{me:witt}. In particular, all the arguments in this section are completely parallel to their respective analogues over a fixed number field developed in \cite{me:witt}. We include them here for the purposes of self-containment and readability. As above, let $F$ be a quadratic form in $N \geq 2$ variables. We first state a result on an effective decomposition of a bilinear space into an orthogonal sum of one-dimensional subspaces, i.e. a version of Siegel's lemma for a bilinear space.

\begin{thm} \label{Siegel:me} Let $Z$ be an $L$-dimensional subspace of $(\qbar^N,F)$, $L<N$. Then there exists a basis $\bx_1,...,\bx_L \in \qbar^N$ for $Z$ such that $F(\bx_i,\bx_j) = 0$ for all $i \neq j$, and
\begin{equation}
\label{Siegel:2}
\prod_{i=1}^L H(\bx_i) \leq 3^{\frac{(L-1)^2(L+2)}{4}} \H(F)^{\frac{L(L+1)}{2}} H(Z)^L.
\end{equation}
\end{thm}

\proof
We argue by induction on $L$. First suppose that $L=1$, then pick any $\boldsymbol 0 \neq \bx_1 \in Z$, and observe that $H(\bx_1) = H(Z)$. Now assume that $L>1$ and the theorem is true for all $1 \leq j < L$. Let $\boldsymbol 0 \neq \bx_1 \in Z$ be a vector guaranteed by Theorem \ref{roy:thunder} so that
\begin{equation}
\label{s2}
H(\bx_1) \leq 3^{\frac{L-1}{2}} H(Z)^{\frac{1}{L}}.
\end{equation}
First assume that $\bx_1$ is a non-singular point in $Z$. Then
$$Z_1 = \{ \bwy \in Z : \bx_1^t F \bwy = 0 \} = \{ \bx_1 \}^{\perp} \cap Z,$$
has dimension $L-1$; here $\{ \bx_1 \}^{\perp} = \{ \bwy \in \qbar^N : \bx_1^t F \bwy = 0 \}$. Then by Lemma \ref{intersection}, Lemma \ref{mult}, and (\ref{s2}) we obtain
\begin{equation}
\label{s3}
H(Z_1) \leq H(\bx_1^t F) H(Z) \leq \H(F) H(\bx_1) H(Z) \leq 3^{\frac{L-1}{2}} \H(F) H(Z)^{\frac{L+1}{L}}.
\end{equation}
Since $\dim_{\qbar}(Z_1) = L-1$, the induction hypothesis implies that there exists a basis $\bx_2,...,\bx_{L}$ for $Z_1$ such that $F(\bx_i,\bx_j) = 0$ for all $2 \leq i \neq j \leq L$, and
\begin{eqnarray}
\label{s4} 
\prod_{i=2}^L H(\bx_i) & \leq & 3^{\frac{(L-2)^2(L+1)}{4}} \H(F)^{\frac{L(L-1)}{2}} H(Z_1)^{L-1} \nonumber \\
& \leq &  3^{\frac{(L-2)^2(L+1)}{4} + \frac{(L-1)^2}{2}} \H(F)^{\frac{L^2+L-2}{2}} H(Z)^{\frac{L^2-1}{L}},
\end{eqnarray}
where the last inequality follows by (\ref{s3}). Combining (\ref{s2}) and (\ref{s4}) we see that $\bx_1,...,\bx_L$ is a basis for $Z$ satisfying (\ref{Siegel:2}) such that $F(\bx_i,\bx_j) = 0$ for all $1 \leq i \neq j \leq L$.

Now assume that $\bx_1$ is a singular point in $Z$. Since $\bx_1 \neq 0$, it must be true that $x_{1j} \neq 0$ for some $1 \leq j \leq N$. Let
$$Z_1 = Z \cap \{ \bx \in \qbar^N : x_j = 0 \},$$
then $\bx_1 \notin Z_1$, $Z = \qbar\bx_1 \perp Z_1$, and 
\begin{equation}
\label{s5}
H(Z_1) \leq H(Z),
\end{equation}
by Lemma \ref{intersection}. Since $\dim_{\qbar}(Z_1) = L-1$, we can apply induction hypothesis to $Z_1$, and proceed the same way as in the non-singular case above. Since the upper bound of (\ref{s5}) is smaller than that of (\ref{s3}), the result follows.
\endproof
\bigskip

Next we discuss the effective structure of the isometry group of a bilinear space over $\qbar$. For the rest of this section assume that $Z \subseteq \qbar^N$ is an $L$-dimensional subspace, $1 \leq L \leq N$ such that the bilinear space $(Z,F)$ is regular. From here on our notation is the same as in section~5 of \cite{me:witt}; we review it here.

First notice that $\qbar^N = Z \perp Z^{\perp_{\qbar^N}}$, where $Z^{\perp_{\qbar^N}} = \{ \bx \in \qbar^N : F(\bx,\bz) = 0\ \forall\ \bz \in Z\}$. Let $\OO(Z,F)$ be the group of isometries of $(Z,F)$, and write $id_Z$ for its identity element. Also let $-id_Z$ be the element of $\OO(Z,F)$ that takes $\bx$ to $-\bx$ for each $\bx \in Z$. Each element $\sigma$ of the isometry group $\OO(\qbar^N,F)$ is uniquely represented by an $N \times N$ matrix $A \in GL_N(\qbar)$, and so we can define $\H(\sigma)$ to be the height of $A$ viewed as a vector in $\qbar^{N^2}$, same way as for the coefficient matrix of~$F$.

Notice that each $\sigma \in \OO(Z,F)$ can be extended to an isometry $\hs \in \OO(\qbar^N,F)$ by selecting an isometry $\sigma' \in \OO(Z^{\perp_{\qbar^N}},F)$. For each $\sigma \in \OO(Z,F)$ choose such an extension $\hs: \qbar^N \rightarrow \qbar^N$ so that $\H(\hs)$ is minimal, and define $\H(\sigma) = \H(\hs)$ for this choice of $\hs$. This definition of height in particular insures that for each $\sigma \in \OO(Z,F)$
\begin{equation}
\label{pm}
\H(\sigma) = \H(-\sigma),
\end{equation}
where $-\sigma = -id_Z \circ \sigma$. Moreover, if $A$ is the matrix of $\hs$, then
\begin{equation}
\label{det}
\det(A) = \det (\hs) = \det \left( \hs \mid_Z \right) \det \left( \hs \mid_{Z^{\perp_{\qbar^N}}} \right) = \det(\sigma) \det(\sigma') = \pm 1.
\end{equation}
We will also refer to this matrix $A$ as the matrix of $\sigma$.

For each $\bx \in Z$ such that $F(\bx) \neq 0$ we can define an element of $\OO(Z,F)$, $\tau_{\bx} : Z \longrightarrow Z$, given by
\begin{equation}
\label{cd1}
\tau_{\bx} (\bwy) = \bwy - \frac{2 F(\bx,\bwy)}{F(\bx)} \bx,
\end{equation}
which is a {\it reflection} in the hyperplane $\{\bx\}^{\perp} = \{ \bz \in Z : F(\bx,\bz) = 0 \}$. It is not difficult to see that the matrix of such a reflection is of the form $(\tau_{ij}(\bx))_{1 \leq i,j \leq N}$, where
\[ \tau_{ij}(\bx) = \left\{ \begin{array}{ll}
    1 - \frac{2}{F(\bx)} \sum_{k=1}^N f_{ik} x_i x_k & \mbox{if $i=j$} \\
    - \frac{2}{F(\bx)} \sum_{k=1}^N f_{jk} x_i x_k & \mbox{if $i \neq j$}
    \end{array}
\right. \]
For each reflection $\tau_{\bx}$, $\det(\tau_{\bx}) = -1$. We say that $\sigma$ is a {\it rotation} if $\det(\sigma)=+1$.
\smallskip

We can now derive some bounds on height of isometries of $(Z,F)$. We start with a simple result, which is precisely Lemma 5.1 of \cite{me:witt}.

\begin{lem} \label{ref_height} Let $\bx \in Z$ be anisotropic and $\tau_{\bx} \in \OO(Z,F)$ be the corresponding reflection. Then
\begin{equation}
\label{ref_bound}
\H(\tau_{\bx}) \leq N^3(N+2) \H(F) H(\bx)^2.
\end{equation}
\end{lem}

An immediate consequence of Lemmas \ref{ref_height} and \ref{nonzero} is the following statement about the existence of a reflection of relatively small height in $\OO(Z,F)$ - this is a direct analogue of Corollary 5.3 of \cite{me:witt}.

\begin{cor} \label{isom:masser} There exists a reflection $\tau \in \OO(Z,F)$ with
\begin{equation}
\label{isom_bound}
\H(\tau) \leq 3^{\frac{(L+2)(L-1)}{2}} 4L N^3 (N+2) \H(F) H(Z)^{\frac{L+2}{L}}.
\end{equation}
\end{cor}

\proof
Let $\bx$ be an anisotropic point in $Z$ guaranteed by Lemma \ref{nonzero}. Let $\tau=\tau_{\bx}$. The result follows by combining (\ref{ref_bound}) with (\ref{anis}).
\endproof

Moreover, every isometry $\sigma \in \OO(Z,F)$ can be represented as a product of reflections of bounded height. This is an effective version of the well-known Cartan-Dieudonn{\'e} theorem. Specifically, we can state the following.

\begin{thm} \label{CD_me} Let $(Z,F)$ be a regular symmetric bilinear space over $\qbar$ with $Z \subseteq \qbar^N$ of dimension $L$, $1 \leq L \leq N$, $N \geq 2$. Let $\sigma \in \OO(Z,F)$. Then either $\sigma$ is the identity, or there exist an integer $1 \leq l \leq 2L-1$ and reflections $\tau_1,...,\tau_l \in \OO(Z,F)$ such that 
\begin{equation}
\label{CD_me_0}
\sigma = \tau_1 \circ \dots \circ \tau_l, 
\end{equation}
and for each $1 \leq i \leq l$,
\begin{equation}
\label{CD_me_1}
\H(\tau_i) \leq \left\{ \left( 2 N^2 3^{\frac{L-1}{2}} \right)^{\frac{L^2}{2}} \H(F)^{\frac{L}{3}} H(Z)^{\frac{L}{2}} \H(\sigma) \right\}^{5^{L-1}}.
\end{equation}
\end{thm}

\noindent
To prove Theorem \ref{CD_me} we will need the following two technical lemmas, which are Lemmas 5.4 and 5.6 of \cite{me:witt}, respectively.

\begin{lem} \label{matrix} Let $A \in GL_N(\qbar)$ be such that $\det(A) = \pm 1$, and write $I_N$ for the $N \times N$ identity matrix. Then
\begin{equation}
\label{matrix_height}
\H(A \pm I_N) \leq 2 \H(A).
\end{equation}
\end{lem}

\begin{lem} \label{matrix1} Let $A$ and $B$ be two $N \times N$ matrices with entries in $\qbar$. Then
\begin{equation}
\label{mht}
\H(AB) \leq \H(A) \H(B).
\end{equation}
\end{lem}
\smallskip

\noindent
{\it Proof of Theorem \ref{CD_me}.}
We argue by induction on $L$. When $L=1$, $Z=\qbar\bx$ for some anisotropic vector $\bx \in \qbar^N$, since $(Z,F)$ is regular. Then $\sigma = \pm id_Z$, where $-id_Z = \tau_{\bx}$, and $\H(\sigma) = \sqrt{N}$ by (\ref{pm}). 
\smallskip

Then assume $L>1$. Write $A$ for the $N \times N$ matrix of $\sigma$, and $I_N$ for the $N \times N$ identity matrix, so in particular $\H(\sigma)=\H(A)$. Notice that for each $\bx \in Z$,
\begin{equation}
\label{zero1}
F(\sigma(\bx)-\bx,\sigma(\bx)+\bx) = 0.
\end{equation}
Let $\bx \in Z$ be the anisotropic vector guaranteed by Lemma \ref{nonzero} with $\sigma(\bx) \pm \bx$ also anisotropic. For this choice of $\pm$, $\tau_{\sigma(\bx) \pm \bx}$ fixes $\sigma(\bx) \mp \bx$ and maps $\sigma(\bx) \pm \bx$ to $-(\sigma(\bx) \pm \bx)$. Then $2 \sigma(\bx) = (\sigma(\bx)) + (\sigma(\bx)-\bx)$ will be mapped to $(\sigma(\bx) \mp \bx)-(\sigma(\bx) \pm \bx) = \mp 2\bx$. We can therefore observe that if $\sigma(\bx) - \bx$ is anisotropic, then
\begin{equation}
\label{sprime1}
\sigma' = \tau_{\sigma(\bx) - \bx} \circ \sigma
\end{equation}
fixes $\bx$. If, on the other hand, $\sigma(\bx) + \bx$ is anisotropic, then
\begin{equation}
\label{sprime2}
\sigma' = \tau_{\sigma(\bx) + \bx} \circ \tau_{\sigma(\bx)} \circ \sigma
\end{equation}
fixes $\bx$. In any case, $\sigma'$ defined either by (\ref{sprime1}) or (\ref{sprime2}) is an isometry of the $(L-1)$-dimensional regular bilinear space $(\{ \bx \}^{\perp},F)$, where $\{ \bx \}^{\perp} = \{ \bz \in Z : F(\bx,\bz) = 0 \}$. Then, by the induction hypothesis,
$$\sigma' = \tau_1 \circ \dots \circ \tau_l,$$
for some reflections $\tau_1,...,\tau_l$ with $1 \leq l \leq 2L-3$ and
\begin{equation}
\label{ind}
\H(\tau_i) \leq \left\{ \left( 2 N^2 3^{\frac{L-2}{2}} \right)^{\frac{(L-1)^2}{2}} \H(F)^{\frac{L-1}{3}} H \left( \{ \bx \}^{\perp} \right)^{\frac{L-1}{2}} \H(\sigma') \right\}^{5^{L-2}},
\end{equation}
for each $1 \leq i \leq l$, and so
\begin{equation}
\label{tau}
\sigma = \sigma'' \circ \tau_1 \circ \dots \circ \tau_l,
\end{equation}
for the same $\tau_1,...,\tau_l$ and $\sigma'' = \tau_{\sigma(\bx) - \bx}$ or $\sigma'' = \tau_{\sigma(\bx) + \bx} \circ \tau_{\sigma(\bx)}$, depending on which of $\sigma(\bx) \pm \bx$ is anisotropic, so $\sigma$ is a product of at most $2L-1$ reflections. Next we are going to produce bounds on their heights. Combining Lemma \ref{ref_height} with a bound analogous to that of Lemma \ref{mult} and with Lemma \ref{nonzero}, we obtain
\begin{equation}
\label{cd6}
\H(\tau_{\sigma(\bx)}) \leq 4L N^3 (N+2)\ 3^{\frac{(L+2)(L-1)}{2}} \H(F) H(Z)^{\frac{L+2}{L}} \H(\sigma)^2.
\end{equation}
Therefore $\tau_{\sigma(\bx)}$ satisfies (\ref{CD_me_1}). Also by Lemma \ref{ref_height},
\begin{equation}
\label{cd10}
\H(\tau_{\sigma(\bx) \pm \bx}) \leq N^3(N+2) \H(F) H(\sigma(\bx) \pm \bx)^2.
\end{equation}
Notice that $\sigma(\bx) \pm \bx = (A \pm I_N) \bx$. Then, once again, by a bound analogous to that of Lemma \ref{mult}
\begin{equation}
\label{cd11}
\H(\sigma(\bx) \pm \bx) \leq H(\bx) \H(A \pm I_N) \leq 2 \sqrt{L}\ 3^{\frac{(L+2)(L-1)}{4}} H(Z)^{\frac{L+2}{2L}} \H(A \pm I_N),
\end{equation}
where the last inequality follows by (\ref{anis}). Combining (\ref{cd11}) with Lemma \ref{matrix}, we obtain
\begin{equation}
\label{cd12}
H(\sigma(\bx) \pm \bx) \leq 4 \sqrt{L}\ 3^{\frac{(L+2)(L-1)}{4}} H(Z)^{\frac{L+2}{2L}} \H(A).
\end{equation}
Combining (\ref{cd10}) and (\ref{cd12}), we obtain
\begin{equation}
\label{cd13}
\H(\tau_{\sigma(\bx) \pm \bx}) \leq 16 L N^3 (N+2)\ 3^{\frac{(L+2)(L-1)}{2}}  \H(F) H(Z)^{\frac{L+2}{L}} \H(\sigma)^2,
\end{equation}
\smallskip
hence $\tau_{\sigma(\bx) \pm \bx}$ satisfies (\ref{CD_me_1}). By combining (\ref{sprime1}), (\ref{sprime2}), (\ref{pm}), Lemma \ref{matrix1}, (\ref{cd6}), and (\ref{cd13}), we have
\begin{equation}
\label{cd13.1}
\H(\sigma') \leq 64 L^2 N^6 (N+2)^2\ 3^{(L+2)(L-1)} \H(F)^2 H(Z)^{\frac{2L+4}{L}} \H(\sigma)^5.
\end{equation}
By Lemma \ref{intersection}, Lemma \ref{mult}, and (\ref{anis})
\begin{equation}
\label{cd13.2}
H \left( \{\bx\}^{\perp} \right) \leq \H(F) H(\bx) H(Z) \leq 2 \sqrt{L}\ 3^{\frac{(L+2)(L-1)}{4}} \H(F) H(Z)^{\frac{3L+2}{2L}}.
\end{equation}
Then bound (\ref{CD_me_1}) follows upon combining (\ref{ind}) with (\ref{cd13.1}) and (\ref{cd13.2}) while keeping in mind that $2 \leq L \leq N$ and $N+2 \leq 2N$. This completes the proof.
\boxed{ }
\bigskip

{\bf Acknowledgment.} I would like to thank Professor Paula Tretkoff for her helpful remarks on the subject of this paper.
\bigskip

\bibliographystyle{plain}  
\bibliography{quad}        

\end{document}